\documentclass{amsart}

\newtheorem{theorem}{Theorem}[section]

\newtheorem{proposition}[theorem]{Proposition}

\newtheorem{_conjecture}[theorem]{Conjecture}

\newtheorem{_problem}[theorem]{Problem}

\newtheorem{_claim}[theorem]{Claim}

\newtheorem{_subclaim}[theorem]{Sub-claim}

\newtheorem{defini}[theorem]{Definition}

\newtheorem{rem}[theorem]{\it Remark}
\newenvironment{remark}{\begin{rem}\rm}{\end{rem}}

\newtheorem{_example}[theorem]{Example}

\newtheorem*{maintheorem}{Main Theorem}

\numberwithin{equation}{section}
\numberwithin{table}{section}
\numberwithin{figure}{section}

\newcommand{\F}{\mathord{\mathbb F}}
\renewcommand{\P}{\mathord{\mathbb  P}}
\newcommand{\Q}{\mathord{\mathbb  Q}}
\newcommand{\R}{\mathord{\mathbb R}}
\newcommand{\Z}{\mathord{\mathbb Z}}

\newcommand{\EEE}{\mathord{\mathcal E}}

\newcommand{\maprightsp}[1]{\; \smash{\mathop{\; \longrightarrow \; }\limits\sp{#1}}\; }

\newcommand{\mapdownsurj}{
\hbox{$\bigm\downarrow$}
\llap{\hbox{\raise 2pt\hbox{$\bigm\downarrow$}}}%
\vstrechmapdown
}

\newcommand{\inj}{\hookrightarrow}

\newcommand{\isom}{\smash{\mathop{\;\to\;}\limits\sp{\sim\,}}}

\newcommand{\set}[2]{\{\; {#1} \; \mid \; {#2} \;  \}}

\newcommand{\map}[3]{ #1 \, : \, #2 \, \to \, #3}

\newcommand{\shortmapisom}[3]{ #1 : #2 \isom #3}

\newcommand{\st}{\subset}

\newcommand{\sprime}{\sp\prime}

\newcommand{\spar}[1]{\sp{(#1)}}

\newcommand{\sperp}{\sp{\perp}}

\newcommand{\inv}{\sp{-1}}

\renewcommand{\qed}{\hfill {$\Box$}}

\newcommand{\pr}{\operatorname{\rm pr}\nolimits}

\newcommand{\Hom}{\operatorname{\rm Hom}\nolimits}

\newcommand{\disc}{\operatorname{\rm disc}\nolimits}

\newcommand{\rmand}{\textrm{and}}

\newcommand{\quand}{\quad\rmand\quad}

\newcommand{\uu}{\mathord{\textbf{\itshape u}}}
\newcommand{\aaa}{\mathord{\textbf{\itshape a}}}
\newcommand{\vv}{\mathord{\textbf{\itshape v}}}
\newcommand{\ww}{\mathord{\textbf{\itshape w}}}

\begin{document}

\title[Supersingular $K3$ surfaces]{
Supersingular $K3$ surfaces in odd characteristic and sextic  double  planes}

\author{Ichiro Shimada}
\address{
Division of Mathematics, 
Graduate School of Science, 
Hokkaido University,
Sapporo 060-0810,
JAPAN
}
\email{shimada@math.sci.hokudai.ac.jp
}


\begin{abstract}
We show that every supersingular $K3$ surface 
is birational to a double cover of a projective plane.
\end{abstract}

\maketitle
\section{Introduction}\label{sec:intro}
A $K3$ surface is called \emph{supersingular}
(in the sense of Shioda)
if the rank of the  Picard lattice  is  equal to $22$.
Supersingular $K3$ surfaces exist only when the characteristic $p$ of the base field is positive.
\par
We  call 
a pair $(X, L)$
of a $K3$ surface $X$ and a line bundle $L$ on $X$
a \emph {sextic double plane}
if $L$ is effective of degree  $L^2=2$, and
$|L|$ has no fixed components.
Let $(X, L)$ be a sextic double plane,
and let
$$
X\;\maprightsp{f}\; Y \;\maprightsp{\pi}\; \P\sp 2
$$
be the Stein factorization of the morphism associated with $|L|$.
Then $f$ is birational, and $\pi$ is a finite morphism of degree $2$.
\par
The purpose of this paper is to prove the following:
\begin{maintheorem}\label{thm:main}
Every supersingular $K3$ surface $X$ has a line bundle $L$
of degree $2$ such that $(X, L)$ is a sextic double plane.
\end{maintheorem}
In~\cite{Shimada2002},
we have proved that every supersingular $K3$ surface in characteristic $2$ is birational to
an \emph{inseparable} double cover of a projective plane.
Therefore we  consider the case where   $p$ 
is an odd prime.
\par
Let $X$ be a supersingular $K3$ surface in characteristic $p$.
We denote by $S\sb X$ the Picard lattice of $X$,
which is  an even lattice of signature $(1,21)$.
Artin~\cite{Artin74} proved  that
the discriminant of  $S\sb X$ 
is equal to $-p\sp{2\sigma\sb X}$,
where $\sigma\sb X$ is a positive integer $\le 10$.
The integer $\sigma\sb X$
is called the \emph{Artin invariant} of $X$.
In~\cite{Artin74, RS_char2,  Shioda79},
it was established that,
for every  pair of a prime integer $p$ and an integer  $\sigma$
with $1\le \sigma\le 10$,
there exists a supersingular $K3$ surface $X$ in characteristic $p$ with $\sigma\sb X=\sigma$.
On the other hand, 
Rudakov and {\v{S}}afarevi{\v{c}}~\cite{Rudakov_Shafarevich81} showed that
the Picard lattice of a supersingular $K3$ surface
is determined, up to isomorphisms of lattices,
by the characteristic $p$ and the Artin invariant.
Moreover, they constructed,
for each $p$ and $\sigma$,  a lattice $\Lambda\sb{p, \sigma}$
that is isomorphic to the Picard lattice of a supersingular $K3$ surface $X$
in characteristic $p$ with $\sigma\sb X=\sigma$.
\par
Let $(X, L)$ be a sextic double plane in characteristic $p$.
When $p> 2$, the double covering  $\pi: Y\to \P\sp 2$ 
induced from $|L|$ is  \emph{separable},
and it branches along 
a plane curve $B\sb{(X, L)}$ of degree $6$
with only rational double points.
Let us denote  by $R\sb{(X, L)}$ the $ADE$-type of the singularities of
$B\sb{(X, L)}$.
\par
In this paper, 
we  give, for each odd prime $p$ and the Artin invariant $\sigma$,
a vector $h\in \Lambda\sb{p, \sigma}$
such that,
under a certain isomorphism
$\shortmapisom{\phi}{\Lambda\sb{p, \sigma}}{S\sb X}$
of lattices,
the line bundle $L$ 
corresponding to  $\phi (h)\in S\sb X$
makes $X$ a sextic double plane.
We also calculate the $ADE$-type $R\sb{(X, L)}$
of some supersingular sextic double planes  $(X, L)$
obtained from the vector $h\in \Lambda\sb{p, \sigma}$.
\par
The proof of Main Theorem depends on elementary but tedious calculations
of linear algebra and quadratic polynomials.
The reader is strongly recommended to use some computer algebra system
while reading this paper.
In particular,
a program that calculates the projection of a given quadratic form
to a coordinate axis (see \S\ref{sec:comp} for the definition)
will be very useful.
\par
\medskip
We will assume from now on that $p$ is an \emph{odd} prime.
\par
\medskip
The author would like to thank Professor I.~R.~Dolgachev
for  stimulating  discussions.
\section{Preliminaries}
We use the following notation and terminologies.
Let $T$ be a lattice;
that is, $T$ is a free $\Z$-module of finite rank
with a non-degenerate symmetric bilinear form
$(x, y)\mapsto xy$
that takes values in $\Z$.
We can express a lattice $T$ by a symmetric matrix $M\sb T$ of  integer components
with respect to a certain basis of $T$.
The \emph{opposite lattice} $T\sp{-}$ of a lattice $T$ is a lattice
with the same underlying $\Z$-module as $T$ such that $M\sb{T\sp -}=-M\sb T$.
The non-degenerate bilinear form induces  a natural injective homomorphism 
$T\inj \Hom (T, \Z).$
The \emph{discriminant group} of $T$ is
the   group $\Hom (T, \Z)/T$,
which is a finite abelian group of order $|\disc T|=|\det M\sb T|$.
For a vector $h\in T$,
we put 
\begin{eqnarray*}
h\sperp &:=& \set{u\in T}{ uh=0}, \\
\EEE\sp{\pm}\sb 0 (h) &:=& \set{u\in T}{uh=\pm 1, u^2= 0}, \\
\EEE\sp{\pm}\sb{\le} (h) &:=& \set{u\in T}{uh=\pm 1, u^2\le 0}.
\end{eqnarray*}
A lattice $T$ is said to be \emph{even}
if $v^2\in 2\Z$ holds for every $v\in T$.
Let $T$ be a (positive or negative) definite even lattice.
A vector $v\in T$ is called a \emph{root} of $T$
if $v^2$ is equal to $2$ or $-2$.
The  roots in $T$ form a root system 
(\cite{Ebeling, Bourbaki}).
We denote by $\Sigma (T)$ the $ADE$-type of this root system,
which in a finite formal sum of the symbols $A\sb l (l\ge 1)$, $D\sb m (m\ge 4)$ and 
$E\sb n (n=6,7,8)$. For simplicity, we put
$$
D\sb 3=A\sb 3, \quad D\sb 2=2A\sb 1, \quand D\sb 1=D\sb 0=0.
$$
See \cite[Chapter 4, Section 7]{Conway_Sloane} for the reason of this convention.
\par
\medskip
Let $X$ be a $K3$ surface in arbitrary characteristic,
and let $S\sb X$ be the Picard lattice of $X$.
\begin{proposition}[\cite{Nikulin90, Urabe86}]\label{prop:urabe}
Let $L$ be a nef line bundle on $X$ with degree $2$.
\par
{\rm (1)}
The pair
$(X, L)$ is a sextic double plane
if and only if the set $\EEE\sp{+}\sb 0([L])$
is empty, where $[L]\in S\sb X$ is the isomorphism class of $L$.
\par
{\rm (2)}
If $(X, L)$ is a sextic double plane,
then $R\sb{(X, L)}$ is equal to $\Sigma ([L]\sperp)$.
\qed
\end{proposition}
\begin{proposition}[\cite{Rudakov_Shafarevich81}]\label{prop:isometry}
Let $v\in S\sb X$ be a vector with $v^2>0$.
Then there exists  an isometry $\gamma : S\sb X\isom S\sb X$
of $S\sb X$ such that $\gamma (v)$ is the isomorphism class of a nef line bundle on $X$.
\qed
\end{proposition}
Let $p$ be an odd prime,
and $\sigma$ a positive integer $\le 10$.
We denote by $\Lambda\sb{p, \sigma}$
a lattice of rank $22$ with the following properties:
\begin{itemize}
\item[(a)]
$\Lambda\sb{p, \sigma}$ is even, 
\item[(b)] the signature of $\Lambda\sb{p, \sigma}$ is  $(1, 21)$, and
\item[(c)]
the discriminant group of $\Lambda\sb{p, \sigma}$
is isomorphic to $(\Z/p\Z) \sp{\oplus 2\sigma}$.
\end{itemize}
\begin{proposition}[\cite{Rudakov_Shafarevich81}]\label{prop:characterization}
These three  conditions determine the lattice $\Lambda\sb{p, \sigma}$
uniquely up to isomorphisms.
\qed
\end{proposition}
\begin{proposition}[\cite{Artin74, Rudakov_Shafarevich81}]\label{prop:isom}
The Picard lattice $S\sb X$ of 
a supersingular $K3$ surface
$X$ in characteristic $p>2$
with  $\sigma\sb X=\sigma$ is isomorphic to 
$\Lambda\sb{p,\sigma}$.
\qed
\end{proposition}
Combining these propositions and 
changing the sign from $\Lambda\sb{p, \sigma}$ to $\Lambda\sb{p, \sigma}\sp -$, 
we obtain the following:
\begin{proposition}\label{prop:lattice}
Let $X$ be a supersingular $K3$ surface
in characteristic $p>2$ with $\sigma\sb X=\sigma$,
and $R$ an $ADE$-type.
There exists a line bundle $L$ on $X$ such that
$(X, L)$ is a sextic double plane
with $R\sb{(X, L)}=R$
if and only if there exists a vector $h\in \Lambda\sb{p, \sigma}\sp-$
such that $h^2=-2$, $\Sigma (h\sperp)=R$,
and $\EEE\sp{-}\sb 0 (h)=\emptyset$.
\qed
\end{proposition}
We  say that \emph{a sextic double plane $(X,  L)$
is obtained from  $h\in \Lambda\sb{p, \sigma}\sp-$}
if $[L]\in S\sb X$ is equal to $\phi\sp - (h)$,
where
$\phi\sp - : \Lambda\sb{p, \sigma}\sp-\isom S\sb X$
is an \emph{anti-}isometry of lattices.
Note that $R\sb{(X, L)}=\Sigma (h\sperp)$ depends not on the choice of $\phi\sp -$
but only on $h\in \Lambda\sb{p, \sigma}\sp-$.
\section{Constructing the Picard lattice}\label{sec:constr}
According to~\cite{Rudakov_Shafarevich81},
we construct the opposite lattice $\Lambda\sb{p, \sigma}\sp-$
of $\Lambda\sb{p, \sigma}$
from the following lattices.
\par
\medskip
(I)
We denote by $U$ and $U\spar{p}$ the even indefinite 
lattices of rank $2$
whose intersection matrices are given by
\begin{equation*}\label{eq:Umat}
\begin{bmatrix}
0 & 1 \\
1 & 0
\end{bmatrix}
\quad
\textrm{and}
\quad
\begin{bmatrix}
0 & p \\
p & 0
\end{bmatrix},
\end{equation*}
with respect a certain basis $\uu\sb 1$ and $\uu\sb 2$, 
respectively.
\par
\medskip
(II)
Let $H\spar{p}$
denote the maximal order in the quaternion algebra
$A\spar{p}$
over $\Q$ ramified only at $p$ and $\infty$.
Then 
$H\spar{p}$
has  a symmetric bilinear form
$$
xy:=\mathord{\textrm{Tr}} (x y\sp*),
$$
where $y\mapsto y\sp*$ is the usual involution of the quaternion
algebra.
With this bilinear form,  $H\spar{p}$ is a positive definite lattice of rank $4$.
We can write the intersection matrix of $H\spar{p}$
by the following result due to Ibukiyama:
\begin{proposition}[\cite{Ibukiyama}]\label{prop:ibukiyama}
Let $p\sb 1, \dots, p\sb r$
be distinct primes.
We put $\varepsilon :=(-1)\sp r$.
Let $q$ be a prime integer 
satisfying $\varepsilon q \equiv 5 \bmod 8$
and
$$
\Bigl(\frac{\varepsilon q }{p\sb i}\Bigr) =-1
\quad \textrm{for all}\quad p\sb i\ne 2.
$$
We put 
$\alpha:=\varepsilon p\sb 1\cdots p\sb r$, $\beta :=\varepsilon q$,
and define a quaternion algebra
$A\spar{p\sb1, \dots, p\sb r}$
to be
$$
\Q + \Q a +\Q b +\Q ab,
\quad\textrm{with}\quad
a^2=\alpha, \;b^2=\beta, \; ab+ba=0.
$$
Then $A\spar{p\sb1, \dots, p\sb r}$
is ramified only at $p\sb 1, \dots, p\sb r$
if $r$ is even,
and only at
$p\sb 1, \dots, p\sb r, \infty$
if $r$ is odd.
Let $\gamma$ be an integer 
satisfying
$\gamma^2 \equiv \alpha \bmod q$.
Then 
\begin{equation*}\label{eq:basisH}
H\spar{p\sb 1, \dots, p\sb r}:=
\Z + \Z \frac{1+b}{2} +\Z \frac{a(1+b)}{2} + \Z \frac{(\gamma +a)b}{q}
\end{equation*}
is the maximal order of $A\spar{p\sb1, \dots, p\sb r}$.
\qed
\end{proposition}

With respect to the basis
$$
\aaa\sb 1:=1, \quad
\aaa\sb 2:=(1+b)/2, \quad
\aaa\sb 3:=a(1+b)/2, \quad
\aaa\sb 4:=(\gamma+a)b/q, 
$$
the intersection matrix of $H\spar{p}$
is written as follows:
\begin{equation*}\label{eq:Hmat}
M\sb {H\spar{p}}:=\begin{bmatrix}
2 & 1 & 0 & 0 \\
1 & {(q+1)}/{2} & 0 & \gamma \\
0 & 0 & {p( q+1)}/{2} & p \\
0 &\gamma & p & {2 (p+\gamma^2) }/{q}
\end{bmatrix},
\end{equation*}
where $q$ is a prime such that 
$$
q \equiv 3 \bmod 8 \quand (\frac{-q}{p})=-1,
$$
and $\gamma$ is an integer  satisfying 
$$
\gamma^2+p \equiv 0\bmod q.
$$
(Note that such a prime $q$ and an integer $\gamma$ always exist.)
Hence the lattice $H\spar{p}$  is even and 
the  discriminant of $H\spar{p}$ is $p^2$.
Moreover, 
since  all the entries of $p M\sb {H\spar{p}} \inv$ are integers,
the discriminant group of $H\spar{p}$
is isomorphic to $(\Z/p\Z)\sp{\oplus 2}$.
\par
\medskip
(III)
For integers $r$ and $s$ satisfying $0\le s\le r$ and $0<r$,
we denote by 
$W\spar{p}\sb{r, s}$
the  lattice of rank $r$
generated by
$\ww\sb 1, \dots, \ww\sb r$ such that
$$
\ww\sb i \ww\sb j =\begin{cases}
p & \text{ if $i=j\le s$,}\\
1 & \text{ if $i=j>  s$,}\\
0 & \text{ if $i\ne j$.}
\end{cases}
$$
We can construct an
even lattice  $V\spar{p}\sb{r, s}$
from $W\spar{p}\sb{r, s}$
by the following method of Venkov~\cite{Venkov78}.
Let $V\sb 0$ be the submodule of
$W\spar{p}\sb{r, s}$
generated by the vectors
$\sum\sb{i=1}\sp r x\sb i \ww\sb i$
with $\sum x\sb i \equiv 0 \bmod 2$,
and let $a$ be the vector
$(1/2) \sum\sb{i=1}\sp{r}\ww\sb i\in W\spar{p}\sb{r, s}\otimes (1/2) \Z$.
We put
$$
V\spar{p}\sb{r, s}:= V\sb 0 \cup (a+V\sb 0) \quad\subset\quad  W\spar{p}\sb{r, s}\otimes  (1/2) \Z.
$$
Then  $V\spar{p}\sb{r, s}$
is a positive definite
even lattice  
with discriminant group isomorphic to  $(\Z/p\Z)\sp{\oplus s}$
if and only if the following holds:
$$
p s + (r-s) \equiv 0 \bmod 8.
$$
\begin{proposition}\label{prop:venkov}
Suppose that  $p s + (r-s) \equiv 0 \bmod 8$ holds.
Then the $ADE$-type
$\Sigma (V\spar{p}\sb{r, s})$ 
of  the roots in 
$V\spar{p}\sb{r, s}$
is equal to
$D\sb{r-s}$ 
except for the cases in Table~\ref{table:venkovexceptions}.
\end{proposition}
\begin{table}
\renewcommand{\arraystretch}{1.4}
\begin{tabular}{|c||c|c|c|c|c|}
\hline
$(p, r, s)$ & $(p, 8, 0) $ & $(7,2,1) $ & $(5, 4,1) $ & $(3,4,2) $ & $(3,6,1) $ \\
\hline
$\Sigma (V\spar{p}\sb{r, s})$  & $E\sb 8$  & $A\sb 1$ & $A\sb 4$ & $2A\sb 2$ & $E\sb 6$ \\
\hline
\end{tabular}
\caption{The exceptional cases of $\Sigma (V\spar{p}\sb{r, s})$}
\label{table:venkovexceptions}
\end{table}
\begin{proof}
We have an  orthogonal decomposition
$$
 W\spar{p}\sb{r, s}=
W\sp{[p]} \oplus W\sp{[1]},
$$
where $W\sp{[p]}$ (resp. $W\sp{[1]}$)
is spanned 
by $\ww\sb 1, \dots, \ww\sb s$
(resp. $\ww\sb{s+1}, \dots, \ww\sb r$).
Suppose that $v\in V\sb 0\subset W\spar{p}\sb{r, s}$
is a root,
and let 
$$
v=v\sp{[p]} + v\sp{[1]}\quad (v\sp{[p]}\in W\sp{[p]}, \;v\sp{[1]}\in W\sp{[1]})
$$
be the decomposition of $v$.
Since the norm of $v\sp{[p]}$ is a multiple of $p>2$,
$v\sp{[p]}$ must be $0$,
and hence $v \in W\sp{[1]} \cap V\sb 0$.
Since $W\sp{[1]} \cap V\sb 0$
is the so-called ``checkerboard lattice" of rank $r-s$,
the $ADE$-type of its roots is equal to $D\sb{r-s}$
(\cite[Chapter 4, Section 7]{Conway_Sloane}).
Suppose that $v\in a+V\sb 0$
is a root.
Then $v$ is written as $\sum b\sb i\ww\sb i/2$,
where $b\sb 1, \dots, b\sb r$ are odd integers
satisfying
$$
p(b\sb 1^2+\dots+b\sb s^2)+(b\sb{s+1}^2+\dots+b\sb r^2)=8.
$$
The left-hand side of this equality is $\ge ps +(r-s)$.
Therefore $ps +(r-s)=8$ holds.
The triples $(p, r, s)$ satisfying this equality 
are just the ones in Table~\ref{table:venkovexceptions}.
Thus, except for these cases, we have proved $\Sigma (V\spar{p}\sb{r, s})=D\sb{r-s}$.
The lattice   $V\spar{p}\sb{8, 0}$ is a positive-definite even unimodular lattice  of rank $8$,
and hence 
is the root lattice of type $E\sb 8$.
The $ADE$-types for the other  triples in Table~\ref{table:venkovexceptions}
can be calculated  directly, for example 
by  the method described in Section~\ref{sec:comp}.
\end{proof}
Let us fix the following vectors in $V\spar{p}\sb{r, s}$:
\begin{eqnarray*}
\vv\sb 1 &:=& \ww\sb 1+\ww\sb 2, \\
\vv\sb 2 &:=& (\ww\sb 1+\cdots + \ww\sb r)/2, \\
\vv\sb j &:=& \ww\sb{j-1} + \ww\sb j \quad (j=3, \dots, r-1), \\
\vv\sb r &:=& 2\ww\sb{r}.
\end{eqnarray*}
When $V\spar{p}\sb{r, s}$ is an even lattice, 
the rank $r$ is even and 
the vectors $\vv\sb1, \dots, \vv\sb r$ form a basis of  $V\spar{p}\sb{r, s}$.
The intersection matrix of  $V\spar{p}\sb{r, s}$
with respect to this basis is easily calculated.
\par
\medskip
Using the characterization  of the lattice $\Lambda\sb{p, \sigma}$
(Proposition~\ref{prop:characterization}),
we see that the following lattices are isomorphic to 
$\Lambda\sb{p,\sigma}\sp-$:
\par
\smallskip
\hbox{
\vbox{\tabskip=14pt
\halign {\strut #\hfil & # \hfil\cr
$U\spar{p} \oplus H\spar{p} \oplus V\spar{p}\sb{16, 2\sigma -4}$
& \text{if	($p\equiv 1\bmod 4$ and  $\sigma >1$) or} \cr
&\text{\phantom{if}  ($p\equiv 3\bmod 4$ and $\sigma \equiv 0\bmod 2$)}; \cr
$U \oplus H\spar{p} \oplus V\spar{p}\sb{16, 2\sigma -2}$
& \text{if ($p\equiv 1\bmod 4$  and $\sigma <10$) or} \cr
& \text{\phantom{if} ($p\equiv 3\bmod 4$ and $\sigma \equiv 1\bmod 2$)}; \cr
$U\spar{p}  \oplus V\spar{p}\sb{20, 2\sigma -2}$
& \text{if   $p\equiv 3\bmod 4$ and $\sigma \equiv 0\bmod 2$}; \cr
$U \oplus V\spar{p}\sb{20, 2\sigma}$
& \text{if  $p\equiv 3\bmod 4$ and $\sigma \equiv 1\bmod 2$}. \cr
}}
}
\section{Computational tools}\label{sec:comp}
Let $S\sb 0$ be a finite set of variables, and 
let $Q(X\sb i|X\sb i\in S\sb 0)$ 
be an inhomogeneous quadratic form of variables  $X\sb i \in S\sb 0$
such that the homogeneous part $Q\sb 2$ of degree $2$
is positive definite.
We denote by $\R\spar{S\sb 0}$ the real affine space 
with the set of coordinates being $S\sb 0$.
We denote by $B\sb Q$ the compact subset
$$
\set{(x\sb i | X\sb i\in S\sb 0 )\in \R\spar{S\sb 0}}{ Q(x\sb i|X\sb i\in S\sb 0)\le 0}
$$
of $\R\spar{S\sb 0}$,
and call it the \emph{quadratic body}
associated with the quadratic form $Q$.
For a non-empty subset $S$ of  $S\sb 0$,
we  denote by
$$
\map{\pr\sb S}{ \R\spar{S\sb 0}}{\R\spar{S}}
$$
the natural projection
to the real affine  space $\R\spar{S}$
with the set of coordinates being $S$.
The image $\pr\sb S(B\sb Q)$
of $B\sb Q$ by $\pr\sb S$
is the quadratic body associated with a new quadratic form
$Q\spar{S}(X\sb j | X\sb j\in S)$
of variables $X\sb j \in S$
constructed by the following procedures.
We solve the system
of linear equations
$$
\frac{\partial Q}{\partial X\sb i}=0 \quad (X\sb i\notin S).
$$
Since $Q\sb 2$ is positive definite,
this system of equations has a solution
of the form
$$
X\sb i =\xi\sb i (X\sb j \;|\; X\sb j\in S)\quad (X\sb i\notin S),
$$
where $\xi\sb i$ is an affine-linear form of variables $X\sb j \in S$.
Then  $Q\spar{S}(X\sb j | X\sb j\in S)$
is  the quadratic form 
obtained from $Q$ by substituting
each $X\sb i \notin S$ with $\xi\sb i$.
\par
We denote the new quadratic form  $Q\spar{S}(X\sb j | X\sb j\in S)$ by
$$
\pi (Q(X\sb i|X\sb i\in S\sb 0), S),
$$
and call it the \emph{projection of $Q$ to $\R\sp{(S)}$}.
When the set $S$ consists of a single element $X\sb \nu$,
$Q\sp{\{X\sb\nu\}}(X\sb \nu)=\pi (Q(X\sb i|X\sb i\in S\sb 0),  \{X\sb \nu\})$ is a quadratic
polynomial of one variable. 
We put
$$
J\pi (Q(X\sb i|X\sb i\in S\sb 0), \{X\sb \nu\}) :=
\set{x\in \R}{Q\sp{\{X\sb\nu\}}(x)\le 0},
$$
which is  just the image of the projection of $B\sb Q$ to the $X\sb \nu$-axis.
\par
Suppose that $S\sb 0=\{X\sb 1, \dots, X\sb n\}$,
and consider  
an inhomogeneous  positive definite quadratic form
$Q(X\sb 1, \dots, X\sb n)$.
Let $A:=\Z\alpha\subset \R$ be the $\Z$-module of rank $1$ generated by 
a non-zero real number $\alpha$.
We can list up all points
$(a\sb 1, \dots, a\sb n)\in A\sp n$
satisfying $Q(a\sb 1, \dots, a\sb n)\le 0$ as follows.
We put
$$
S\sb k:=\{X\sb 1, \dots, X\sb k\},
\quad\textrm{and}\quad
Q\sp{(S\sb k)}(X\sb 1, \dots, X\sb k) :=\pi (Q (X\sb 1, \dots, X\sb n), S\sb k).
$$
If we have a point $(a\sb 1, \dots, a\sb k)\in A\sp k$
such that 
$Q\sp{(S\sb k)}(a\sb 1, \dots, a\sb k) \le 0$,
then it is easy to calculate  the set
\begin{multline*}
\phantom{aaaa}\set{a\in A}{Q\sp{(S\sb {k+1})}(a\sb 1, \dots, a\sb k, a)\le 0}=\\
J\pi (Q(a\sb 1, \dots, a\sb k, X\sb{k+1}, \dots, X\sb n), \{X\sb{k+1}\})\cap A,\phantom{aaaa}
\end{multline*}
where
$Q(a\sb 1, \dots, a\sb k, X\sb{k+1}, \dots, X\sb n)$
is considered as a quadratic form of variables $X\sb{k+1}, \dots, X\sb n$.
Starting from $Q\sp{(S\sb 1)}(X\sb 1)$,
we can make inductively  the list of all  points in $B\sb Q \cap  A\sp n$.
\par
Let $T$ be a positive definite even lattice of rank $n$,
and let $Q\sb T(X\sb 1, \dots, X\sb n)$
be a homogeneous quadratic form associated with  $T$.
We can calculate the set of roots of $T$
by applying the above algorithm  to
$Q\sb T(X\sb 1, \dots, X\sb n)-2$.
\section{Proof of Main Theorem}\label{sec:proof}
The strategy of the proof is as follows.
It is enough to find a vector $h\in \Lambda\sb{p, \sigma}\sp-$
such that $h^2=-2$ and $\EEE\sp{-}\sb 0(h)=\emptyset$.
We decompose $\Lambda\sb{p, \sigma}\sp-$
into an orthogonal direct sum
$M\oplus V$,
where
$M$ is of rank $r$ with signature $(r-1, 1)$,
and show that there exists a vector $h\sb 0 \in M$ satisfying $h\sb 0^2=-2$ and 
$\EEE\sp{-}\sb{\le} (h\sb 0)=\emptyset$.
Since $V$ is positive definite,
such a vector $h\sb 0\in M$
 yields the hoped-for vector $h\in \Lambda\sb{p, \sigma}\sp-$
by the natural inclusion $M\inj \Lambda\sb{p, \sigma}\sp-$.

In each of the cases below,
we explicitly give a vector $h\sb 0\in M$
satisfying $h\sb 0^2=-2$,
and a basis $e\sb 1, \dots, e\sb{r-1}\in M$
of $h\sb 0\sperp$.
We put
$$
v\sb 0:=h\sb 0/2\;\;\in\;\; M\otimes  (1/2) \Z,
$$
and define an inhomogeneous quadratic form $Q$ 
with variables $X\sb 1, \dots, X\sb{r-1}$ by 
$$
Q(X\sb 1, \dots, X\sb {r-1}):=(v\sb 0 + X\sb 1 e\sb 1+\dots +X\sb {r-1} e\sb {r-1})^2.
$$
By the \emph{$\Z$-condition},
we mean a necessary and sufficient condition
on  $x\sb 1, \dots, x\sb{r-1} \in  (1/2) \Z$
for a vector $v\sb 0 + x\sb 1 e\sb 1+ \cdots + x\sb {r-1} e\sb {r-1}\in M\otimes  (1/2) \Z$
to be in $M$.
Since $v\sb 0 h\sb 0=-1$,
every vector $x\in M$
satisfying $xh\sb 0=-1$ is uniquely written as
$$
x=v\sb 0 + x\sb 1 e\sb 1+\cdots +x\sb{r-1} e\sb{r-1},
$$
where $x\sb 1, \dots, x\sb{r-1}\in  (1/2) \Z$ satisfy the $\Z$-condition.
\par
In order to prove 
that  $\EEE\sp{-}\sb \le(h\sb 0)$ is empty,
it is therefore enough to show that there are no $x\sb 1, \dots, x\sb{r-1}\in  (1/2) \Z$
satisfying the $\Z$-condition
and the inequality $Q(x\sb1, \dots, x\sb{r-1})\le 0$.
%

We  also investigate the $ADE$-type
$R\sb{(X, L)}$ of the sextic double plane
$(X, L)$
obtained from  the vector $h\in \Lambda\sb{p, \sigma}\sp-$
corresponding to the given $h\sb 0\in M$ via $M\inj \Lambda\sb{p, \sigma}\sp-$.
Since
$$
R\sb{(X, L)}=\Sigma (h\sperp)=\Sigma (h\sb 0\sperp) + \Sigma (V),
$$
it is enough to calculate $\Sigma (V)$
and the set 
of the roots 
in the positive definite even lattice $h\sb 0\sperp$.
In order to calculate this set, 
we define the quadratic form $G$ with variables
$X\sb 1, \dots,  X\sb {r-1}$
by
$$
G(X\sb 1, \dots,  X\sb {r-1}):=-2+ (X\sb 1 e\sb 1+\dots +X\sb {r-1} e\sb {r-1})^2,
$$
and find all the  points $(b\sb 1, \dots, b\sb{r-1})\in \Z\sp{\oplus{(r-1)}}$
satisfying $G(b\sb 1, \dots, b\sb{r-1})=0$.
\par
\medskip
We divide the entire situation into the following overlapping cases:
\par
\smallskip
{\bf Case I.}
$p\equiv 1 \bmod 4$ and $\sigma <10$.
\par
\smallskip
{\bf Case II.}
$p\equiv 1 \bmod 4$ and $\sigma >1$.
\par
\smallskip
{\bf Case III.}
$p\equiv 3 \bmod 4$ and $\sigma\equiv 0 \bmod 2$.
\par
\smallskip
{\bf Case IV.}
$p\equiv 3 \bmod 4$ and $\sigma\equiv 1 \bmod 2$.
%
%
%
%
\subsection{ Case I}
$p\equiv 1 \bmod 4$ and $\sigma <10$.
In this case, we have
$$
\Lambda\sb{p, \sigma}\sp- \;\cong\; U \oplus H\spar{p} \oplus V\spar{p}\sb{16, 2\sigma-2}.
$$
We choose  $U\oplus H\spar{p}$ as the lattice $M$.
We  express  vectors of $U\oplus H\spar{p}$
as  row vectors
with respect to the basis $\uu\sb 1, \uu\sb 2, \aaa\sb1, \dots, \aaa\sb 4$.
Replacing $\gamma$ in the construction of $H\spar{p}$ by $\gamma+q$ if necessary,
we can assume that $\gamma$ is odd.
Then $p+\gamma^2 \equiv 2 \bmod 4$, because $p\equiv 1 \bmod 4$.
Therefore
$$
t:=-\frac{1}{4}\bigl(\frac{p+\gamma^2}{q}+2\bigr)
$$
is an integer.
We put
$$
\hbox{
\def\hs{\hskip 10pt}
\vbox{\tabskip=0pt
\halign {\strut $#$ \hs $:= [$ &\hfil \hs $#$\hs \hfil&\hfil\hs $#$ \hs \hfil&\hfil $#$ \hfil
&\hfil\hs  $#$\hs \hfil &\hfil\hs $#$\hs \hfil &\hfil\hs  $#$\hs  \hfil&\hfil#$]$&# \hfil\cr 
h\sb 0 & 2, & 2\;t, & 1, & 0, & 0, & 1&&, \cr
e\sb 1 & 1, & 0, & -t, & 0, & 0, & 0 &&,\cr 
e\sb 2 & 0, & 1, & -1, & 0, & 0, & 0 &&, \cr
e\sb 3 & 0, & 0, & -(\gamma+1)/2, & 1, & 0, & 0 &&, \cr
e\sb 4 & 0, & 0, & -(p+\gamma^2)/q, & 0, & 0, & 1 &&, \cr
e\sb 5 & 0, & 0, & -p, & 0, & 2, & 0& &.\cr
}}}
$$
It is easy to see that  $h\sb 0\sp 2=-2$, and 
that $e\sb1, \dots, e\sb 5$ form  a basis of $h\sb 0\sperp$,
because $p$ and $2$ are prime to each other.
The $\Z$-condition in this case is as follows:
\begin{equation}\label{eq:ZI}
x\sb 1, x\sb 2, x\sb 3 \in \Z, \quad  x\sb 4, x\sb 5 \in \Z+1/2,
\end{equation}
because $-(p+\gamma^2)/q$ is an even integer.
We assume  that there exist $a\sb 1, \dots, a\sb 5\in  (1/2) \Z$ 
satisfying~\eqref{eq:ZI} and $Q(a\sb 1, \dots, a\sb 5)\le 0$,
and derive a contradiction.
Since
$$
J\pi (Q(X\sb 1, X\sb 2, X\sb 3, X\sb 4, X\sb 5), \{X\sb 1\})=[-1, 1],
$$
we have $a\sb 1=0$ or $a\sb 1=\pm 1$.
From 
\begin{eqnarray*}
J\pi (Q(\pm 1, X\sb 2, X\sb 3, X\sb 4, X\sb 5), \{X\sb 5\})&=&[0, 0],\\
J\pi (Q(\phantom{\pm } \hskip -3pt 0 \hskip 3pt, X\sb 2, X\sb 3, X\sb 4, X\sb 5), \{X\sb
5\})&=&[-\sqrt{p}/2p,
\sqrt{p}/2p],
\end{eqnarray*}
and $a\sb 5\in \Z+1/2$,
we get a contradiction.
Thus $\EEE\sp -\sb{\le} (h\sb 0)=\emptyset$ is proved.
\qed
\par
\medskip
Let us investigate the $ADE$-type 
$$
R\sb{(X, L)}=\Sigma (h\sb 0\sperp) + \Sigma (V\spar{p}\sb{16, 2\sigma -2})=
\Sigma (h\sb 0\sperp)+D\sb{18-2\sigma}
$$
of the sextic double plane $(X, L)$
constructed from $h\sb 0$.
There are at least two roots 
$$
[0,1,-1,0,0,0],
\quand 
[2, \,-{(p+\gamma^2)}/{2q}\, , 1,0,0,1 ]
$$
in $h\sb 0\sperp$,
which are perpendicular to each other.
Numerical experiments show that the $ADE$-type $\Sigma (h\sb 0\sperp)$ depends on the
choice of  $q$ and $\gamma$ in the construction of $H\spar{p}$.
See Table~\ref{table:examplesCaseI}.
\begin{table}
$$
\hbox{
\vbox{\tabskip=0pt \offinterlineskip
\def\tableskip{& height 2pt &\omit& height 2pt &\omit& height 2pt &\omit& height 2pt \cr}
\def\tablerule{\noalign{\hrule}}
\halign {
 #& \vrule#\tabskip=1em plus 2em  & \hfil # \hfil &\vrule#&   \hfil # \hfil &\vrule#&  \hfil # \hfil &\vrule#\tabskip 0pt \cr
\tablerule
\tableskip
&& $p$ && $(q, \gamma)$ && $R\sb{(X, L)}$ & \cr \tableskip\tablerule\tableskip \tablerule\tableskip 
&& $41$ && $(3, 1)$ &&  $A\sb 1+ A\sb 2  + D\sb{18-2\sigma}$ & \cr \tableskip
&&   && $(11, 5)$ && $2 A\sb 1+   D\sb{18-2\sigma}$ & \cr\tableskip\tablerule \tableskip
&& $53$ && $(3, 1)$ &&  $A\sb 1+ A\sb 2  + D\sb{18-2\sigma}$ & \cr \tableskip
&&   && $(67, 9)$ &&  $A\sb 3  + D\sb{18-2\sigma}$ & \cr \tableskip\tablerule \tableskip
&& $61$ && $(11, 7)$ &&  $2A\sb 1  + D\sb{18-2\sigma}$ & \cr \tableskip
&&   && $(43, 5)$ &&  $A\sb 3  + D\sb{18-2\sigma}$ & \cr \tableskip\tablerule\tableskip 
&& $101$ && $(3,1)$ &&  $A\sb 1  + A\sb 2 + D\sb{18-2\sigma}$ & \cr \tableskip
&&   && $(11, 3)$ &&  $2 A\sb 1  + D\sb{18-2\sigma}$ & \cr \tableskip
&&   && $(163, 15)$ &&  $A\sb 3  + D\sb{18-2\sigma}$ & \cr \tableskip\tablerule 
}}
\hfill
}
$$
\caption{Examples of $R\sb{(X, L)}$ in Case I}
\label{table:examplesCaseI}
\end{table}
\begin{remark}
When $p=5$ and $(q, \gamma)=(3,1)$,
there exist  two other   roots 
$$
[0, -1, 0, 1, 0, 0], \quand 
[-1, 0, 0, -1, 0, 0]
$$
in $h\sb 0\sperp$,
and  $\Sigma (h\sb 0\sperp)=A\sb 4$ holds.
Using the isomorphisms of $\Lambda\sb{5, \sigma}\sp-$
with the lattices in Table~\ref{table:exampleschar5},
we obtain examples of sextic curves $B\sb{(X,L)}$
with only rational double points
such that 
 the total Milnor number is $20$.
Note that, in characteristic $0$,  the  maximum 
of the total Milnor number 
of sextic curves with only rational double points   
is $19$.
See Yang~\cite{Yang}.
\begin{table}
$$
\hbox
{
\renewcommand{\arraystretch}{1.2}
\begin{tabular}{|ccc|l|l|}
\hline
&$\sigma$ && \hfil $\Lambda\sp-\sb{5, \sigma}$ & $R\sb{(X, L)}$ \\
\hline
&$1$&  & $U\oplus H\spar{5} \oplus V\spar{5}\sb{16, 0}$\phantom{\vrule height
14pt}  &
$A\sb 4 + D\sb{16}$
\\
&$1$& & $U\oplus H\spar{5} \oplus E\sb 8 \oplus E\sb 8$   & $A\sb 4 + 2 E\sb 8$ \\
&$2$& & $U\oplus H\spar{5}  \oplus A\sb 4 \oplus A\sb 4 \oplus E\sb 8$  & $3A\sb 4
+  E\sb 8$\\
&$3$& & $U\oplus H\spar{5}  \oplus A\sb 4 \oplus A\sb 4 \oplus  A\sb 4 \oplus A\sb
4$ & $5A\sb 4$\\
\hline
\end{tabular}
}
$$
\medskip
\caption{Examples of $B\sb{(X, L)}$
with the total Milnor number $20$  in characteristic $5$}
\label{table:exampleschar5}
\end{table}
\par
On the other hand, there exists a supersingular sextic double plane $(X, L)$
in characteristic $5$
with $\sigma\sb X=1$ such that the branch curve $B\sb{(X, L)}$
is smooth.
Indeed,
let $B\st \P\sp 2$ the Fermat curve
$$
x\sb 0\sp 6 + x\sb 1 \sp 6 + x\sb 2 \sp 6=0
$$
of degree $6$,
and let $X\sb B$ be the double cover of $\P\sp 2$ that branches along $B$.
If  $P\in B$  is an $\F\sb{25}$-rational point of $B$,
then the tangent line $\ell\sb P$ to $B$ at $P$ intersects $B$
only at $P$ with multiplicity $6$,
and hence the pull-back of $\ell\sb P$ to $X\sb B$ splits into
two $(-2)$-curves $\ell\sb P\sp +$ and $\ell\sb P\sp -$
intersecting only at one point with multiplicity $3$.
The number of $\F\sb{25}$-rational points of $B$ is
$126$,
and hence we obtain $252$ smooth rational curves $\ell\sb P\sp\pm$ on $X\sb B$.
It is easy to make the matrix of intersection numbers  between  these curves.
Choosing suitable $22$ curves from them,
we obtain a matrix of determinant $-25$.
Hence $X\sb B$ is supersingular, and the Artin invariant of $X\sb B$ is $1$.
\end{remark}
\subsection{ Case II}
$p\equiv 1 \bmod 4$ and $\sigma >1$.
In this case, we have
$$
\Lambda\sb{p, \sigma}\sp- \;\cong\; U\spar{p} \oplus H\spar{p} \oplus V\spar{p}\sb{16, 2\sigma-4}.
$$
We show that there exists a vector $h\sb 0\in M:=U\spar{p}\oplus H\spar{p}$
such that $h\sb 0\sp 2 =-2$ and $\EEE\sp{-}\sb{\le} (h\sb 0)=\emptyset$.
 Vectors of $U\spar{p}\oplus H\spar{p}$ are expressed as  row vectors 
with respect to the basis $\uu\sb 1, \uu\sb 2, \aaa\sb1, \dots, \aaa\sb 4$.
Since $p\equiv 1 \bmod 4$,
there exists an integer $\alpha$ such that $\alpha^2\equiv -1\bmod p$.
Replacing $\alpha$ with $\alpha+p$
if necessary,
we can assume that $\alpha$ is even.
We put 
$$
b:=-\frac{\alpha^2+1}{p},
$$
 and set
$$
\hbox{
\def\hs{\hskip 10pt}
\vbox{\tabskip=0pt
\halign {\strut $#$ \hs $:= [$ &\hfil \hs $#$\hs \hfil&\hfil\hs $#$ \hs \hfil&\hfil $#$ \hfil
&\hfil\hs  $#$\hs \hfil &\hfil\hs $#$\hs \hfil &\hfil\hs  $#$\hs  \hfil&\hfil#$]$&# \hfil\cr 
h\sb 0 & 1, & b, & \alpha, & 0, & 0, & 0&&, \cr
e\sb 1 & 0, & 0, & 0, & 0, & 0, & 1 &&,\cr 
e\sb 2 & 0, & 0, & 0, & 0, & 1, & 0 &&, \cr
e\sb 3 & 1, & -b, & 0, & 0, & 0, & 0 &&, \cr
e\sb 4 & 0, & 0, & 1, & -2, & 0, & 0 &&, \cr
e\sb 5 & 0, & -\alpha, & 0, & p, & 0, & 0 &&.\cr
}}}
$$
Then $h\sb 0^2$ is equal to $-2$,
and $e\sb 1, \dots, e\sb 5$ form a basis of $h\sb 0\sperp$.
The $\Z$-condition in this case is as follows:
\begin{equation}\label{eq:ZII}
x\sb 1, x\sb 2, x\sb 4, x\sb 5\in \Z
\quand
x\sb 3 \in \Z+1/2,
\end{equation}
because $\alpha$ is even.
We assume  that there exist $a\sb 1, \dots, a\sb 5 \in  (1/2) \Z$ 
satisfying~\eqref{eq:ZII}
and 
 $Q(a\sb 1, \dots, a\sb 5)\le 0$,
and derive a contradiction.
Since
\begin{eqnarray*}
J\pi (Q(X\sb 1, X\sb 2, X\sb 3, X\sb 4, X\sb 5), \{X\sb 2\}) &= &[-1/\sqrt{p}, 1/\sqrt{p}],\\
J\pi (Q(X\sb 1, X\sb 2, X\sb 3, X\sb 4, X\sb 5), \{X\sb 3\}) &= &[-1/2, 1/2],
\end{eqnarray*}
we have $a\sb 2=0$ and $a\sb 3=\pm 1/2$.
Since
$$
\pi (Q(X\sb 1, 0, \pm 1/2, X\sb 4, X\sb 5), \{X\sb 5\})=(p X\sb 5\mp \alpha)^2/2,
$$
we have $a\sb 5=\pm \alpha/p$,
which contradicts $a\sb 5\in \Z$.
Thus $\EEE\sp{-}\sb{\le} (h\sb 0)=\emptyset$ is proved.
\qed
\par
\medskip
We show that $\Sigma (h\sb 0\sperp)=0$ in this case, so that the $ADE$-type 
$R\sb{(X, L)}$
of the sextic double plane $(X, L)$
constructed from $h\sb 0$ above is equal to 
$$
R\sb{(X, L)}=\Sigma (h\sb 0\sperp) + \Sigma (V\spar{p}\sb{16, 2\sigma -4})=
D\sb{20-2\sigma}.
$$
Recall that  
$$
G(X\sb 1, X\sb 2, X\sb 3, X\sb 4 , X\sb 5)=
-2+(X\sb 1 e\sb 1 +X\sb 2 e\sb 2 +X\sb 3 e\sb 3 + X\sb 4 e\sb 4 + X\sb 5 e\sb 5)^2.
$$
We assume that there exist  integers  $b\sb 1, \dots, b\sb 5$
such that $G(b\sb 1, \dots, b\sb 5)=0$, and derive a contradiction.
From
\begin{eqnarray*}
J\pi (G(X\sb 1, X\sb 2, X\sb 3, X\sb 4, X\sb 5), \{X\sb 2\}) &=& [-2/\sqrt{p}, 2/\sqrt{p}],\\
J\pi (G(X\sb 1, X\sb 2, X\sb 3, X\sb 4, X\sb 5), \{X\sb 3\}) &=& [-1, 1],
\end{eqnarray*}
we obtain $b\sb 2=0$ because $p>4$, and $b\sb 3=0$ or $\pm1$.
Suppose that $b\sb 3=0$.
Since
$$
J\pi (G(X\sb 1, 0, 0, X\sb 4, X\sb 5), \{X\sb 5\}) =[-2/p, 2/p], 
$$
we have $b\sb 5$=0.
Then 
$$
G(X\sb 1, 0, 0, X\sb 4, 0)=\frac{2}{q}(
(q X\sb 4 -\gamma X\sb 1)^2 + p X\sb 1 ^2-q
).
$$
Therefore we must have
$(qb\sb 4 -\gamma b\sb 1)^2\equiv q\bmod p$,
which is impossible because $q$ is in  a non-quadratic residue modulo $p$
in the case $p\equiv 1\bmod 4$.
Suppose that $b\sb 3=\pm 1$.
Since
$$
\pi (G(X\sb 1, 0, \pm 1, X\sb 4, X\sb 5), \{X\sb 5\})=(p X\sb 5\mp 2\alpha)^2/2,
$$
we obtain a contradiction to $b\sb 5 \in \Z$.
Therefore the assertion $R\sb{(X, L)}=D\sb{20-2\sigma}$ is proved.
\qed
%
%
%
%
%
%
\subsection{  Case III}
$p\equiv 3 \bmod 4$ and $\sigma\equiv 0 \bmod 2$.
In this case, we have
$$
\Lambda\sb{p, \sigma}\sp- \;\cong\; U\spar{p} \oplus H\spar{p} \oplus V\spar{p}\sb{16, 2\sigma-4}.
$$
We choose $ U\spar{p} \oplus H\spar{p}$ as $M$.
By replacing $\gamma$ in the construction of $H\spar{p}$
with $\gamma+q$ if necessary,
we can assume that
$\gamma\not\equiv 0 \bmod p$.
Using Chevalley-Warning theorem (\cite{Serre_CA})
and the assumption $p\equiv 3 \bmod 4$,
we see that there exists a solution
of the equation 
$$
X^2+\alpha Y^2 =-1
$$
with $Y\ne 0$
in $\F\sb p$ for any $\alpha\in \F\sb p$.
Therefore we have integers $x$ and $y$ such that
$$
x^2 + xy +\frac{q+1}{4}y^2 =\Bigl( x+ \frac{y}{2}\Bigr)^2 +\frac{q}{4} y^2 \equiv -1 \bmod p 
$$
and $y\not\equiv 0\bmod p$.
Replacing $x$ and $y$ with $x+p$ and $y+p$ if necessary,
we can assume that both of 
$x$ and $y$ are even.
We then put
$$
b:=-\frac{1}{p}(x^2 + xy +\frac{q+1}{4}y^2+1)\in \Z.
$$
Note that $b$ is odd.
Since $\gamma y\not \equiv 0\bmod p$,
there exist integers 
$E\sb{36}$ and $E\sb{46}$   satisfying
$$
\gamma y E\sb{36} + (2x +y) \equiv 0\bmod p,
\quad 
\gamma y E\sb{46} + x+ \frac{1+q}{2}y \equiv 0\bmod p.
$$
We then put
$$
E\sb{32}:= -\frac{1}{p}(\gamma y E\sb{36} + 2x + y  )\in \Z,
\quad
E\sb{42}:= -\frac{1}{p}(\gamma y E\sb{46} + x+ \frac{1+q}{2}y  )\in \Z,
$$
and set
$$
\hbox{
\def\hs{\hskip 10pt}
\vbox{\tabskip=0pt
\halign {\strut $#$ \hs $:= [$ &\hfil \hs $#$\hs \hfil&\hfil\hs $#$ \hs \hfil&\hfil $#$ \hfil
&\hfil\hs  $#$\hs \hfil &\hfil\hs $#$\hs \hfil &\hfil\hs  $#$\hs  \hfil&\hfil#$]$&# \hfil\cr 
h\sb 0 & 1, & b, & x, & y, & 0, & 0&&, \cr
e\sb 1 & 0, & 0, & 0, & 0, & 1, & 0 &&,\cr 
e\sb 2 & 1, & -b, & 0, & 0, & 0, & 0 &&, \cr
e\sb 3 & 0, & E\sb{32}, & 1, & 0, & 0, & E\sb{36} &&, \cr
e\sb 4 & 0, & E\sb{42}, & 0, & 1, & 0, & E\sb{46} &&, \cr
e\sb 5 & 0, & -\gamma y, & 0, & 0, & 0, & p&&.\cr
}}}
$$
It is easy to see that
 $h\sb 0\sp 2=-2$,
and that $e\sb1, \dots, e\sb 5$ form  a basis of $h\sb 0\sperp$,
because $p$ and $-\gamma y$ are prime to each other.
The $\Z$-condition in this case is as follows:
\begin{equation}\label{eq:ZIII}
x\sb 1,  x\sb 3, x\sb 4, x\sb 5 \in \Z, \quad   x\sb 2\in \Z+1/2,
\end{equation}
because $x$ and $y$ are even, and $b$ is odd.
Suppose that there exist $a\sb 1, \dots, a\sb 5\in  (1/2) \Z$
satisfying~\eqref{eq:ZIII} and  $Q(a\sb1, \dots, a\sb5)\le 0$.
Since
$$
J\pi ( Q(X\sb1,  X\sb 2, X\sb 3, X\sb 4,  X\sb 5), \{X\sb 1\}) = [-1/\sqrt{p}, 1/\sqrt{p}],
$$
we have $a\sb 1=0$.
Since 
$$
J\pi ( Q( 0, X\sb 2, X\sb 3, X\sb 4,  X\sb 5), \{X\sb 2\}) = [-1/2,1/2],
$$
we have $a\sb 2=\pm 1/2$.
Because
$\pi ( Q(0 , \pm 1/2, X\sb 3, X\sb 4, X\sb 5), \{X\sb 5\})$
is a  multiple of $(x E\sb{36}  + y E\sb{46}  \pm 2pX\sb 5)^2$
by a positive constant, 
we obtain
$$
a\sb 5=\mp\frac{x E\sb{36}  + y E\sb{46} }{2p}.
$$
In $\F\sb p$, however, we have
$$
-(x E\sb{36}  + y E\sb{46} )=\frac{4 x^2 + 4 x y + (1+q) y^2}{2\gamma y}=\frac{-2}{\gamma y}\ne 0.
$$
Therefore $a\sb 5$ cannot be an integer,
and  we arrive at  a contradiction.
Thus $\EEE\sp-\sb{\le}(h\sb 0)=\emptyset$ is proved.
\qed
\par
\medskip
We  prove that
$$
R\sb{(X, L)}=\Sigma (V\spar{p}\sb{16, 2\sigma -4})= D\sb{20-2\sigma}
$$
holds for a sextic double plane $(X, L)$
obtained from the vector $h\sb 0$ above.
It is enough to show that $\Sigma (h\sb 0\sperp)=0$.
Again we assume that there exists an integer point $(b\sb 1, \dots, b\sb 5)\in \Z\sp{\oplus 5}$
satisfying $G(b\sb 1, \dots, b\sb 5)=0$,
and derive a contradiction.
Since
$$
J\pi (G(X\sb1, X\sb 2, X\sb 3, X\sb 4, X\sb 5), \{ X\sb 1\})=[-2/\sqrt{p}, 2/\sqrt{p}],
$$
we see that $b\sb 1=0$ if $p>3$,  and $b\sb 1=0$ or $\pm1$ if $p=3$.
Since 
$$
J\pi (G(X\sb1, X\sb 2, X\sb 3, X\sb 4, X\sb 5), \{ X\sb 2\})=[-1, 1],
$$
we see that $b\sb 2=0$ or $\pm1$.
Suppose that $(b\sb 1, b\sb 2)=(0,0)$.
Then we have
$$
2qy^2 G(0,0,X\sb 3, X\sb 4, X\sb 5)\equiv -4((2X\sb 3+X\sb 4)^2 + y^2 q)\;\;\bmod p.
$$
Since $-q$ is in a non-quadratic residue modulo $p$ and $y\not\equiv 0 \bmod p$,
we get a contradiction.
Suppose that $(b\sb 1, b\sb 2)=(0, \pm 1)$.
The projections of  $G(0, \pm1, X\sb 3, X\sb 4, X\sb 5)$ to 
the $X\sb 3$-axis and the $X\sb4$-axis
are  multiples of $(X\sb 3\mp x)^2$ and
$(X\sb 4\mp y)^2$ by positive constants, respectively.
Therefore 
$(b\sb 3, b\sb 4)=\pm (x, y)$.
Solving the equation $G(0, \pm 1, \pm x, \pm y, X\sb 5)=0$,
we obtain 
$$
b\sb 5=\mp (x E\sb{36}  + y E\sb{46} )/p,
$$
which contradicts  $b\sb 5\in \Z$.
Suppose that $p=3$ and $b\sb 1=\pm1$.
Since
$$
J\pi (G(\pm 1, X\sb 2, X\sb 3, X\sb 4, X\sb 5), \{ X\sb 2\})=[-1/2, 1/2],
$$
we obtain $b\sb 2=0$.
Let $\nu$ be the positive integer $(p+\gamma^2)/q$.
Projecting the quadratic form $G(\pm 1, 0, X\sb 3, X\sb 4, X\sb 5)$ to the $X\sb 3$-axis and the $X\sb
4$-axis, we see that 
$$
3(2b\sb 3 \pm\gamma)^2-\nu-3=0 \quand 3(b\sb 4\mp\gamma)^2-\nu=0.
$$
In particular, both of $\nu/3$ and $\nu/3+1$ are square integers.
Thus we get a contradiction, and the proof of the assertion $R\sb{(X, L)}= D\sb{20-2\sigma}$
is completed.
\qed
\subsection{  Case IV}
$p\equiv 3 \bmod 4$ and $\sigma\equiv 1 \bmod 2$.
We have
$$
\Lambda\sb{p, \sigma}\sp- \;\cong\; U \oplus V\spar{p}\sb{20, 2\sigma}.
$$
In this case, $M$ is the entire lattice $U \oplus V\spar{p}\sb{20, 2\sigma}$.
We  express  vectors of $U \oplus V\spar{p}\sb{20, 2\sigma}$
as row vectors
with respect to the basis $\uu\sb 1, \uu\sb 2, \vv\sb 1, \dots, \vv\sb{20}$,
where
$\vv\sb 1, \dots, \vv\sb{20}$
are the basis of $V\spar{p}\sb{20, 2\sigma}$ fixed  in \S\ref{sec:constr}~(III).
We put
$$
\hbox{
\def\hs{\hskip 5pt}
\vbox{\tabskip=0pt
\halign {\strut $#$ \hs $:= [$ &\hfil \hs $#$\hs \hfil&\hfil\hs $#$ \hs \hfil&\hfil $#$ \hfil
&\hfil\hs  $#$\hs \hfil &\hfil\hs $#$\hs \hfil &\hfil\hs  $#$\hs  \hfil&
\hfil \dots $#$\dots,  &\hfil $#$ \hfil &$]$#& # \cr
h\sb 0 & 2, & -(p+1)/2, & 1,   & 0, & 0, & 0, &          & 0 &&, \cr
e\sb 1 & 0, & 0, 							& 0, 		& 0, & 0, & 0, &										&	1 &&,\cr 
\multispan{11}{\dots\dots\dots\phantom{\vrule height 10pt depth 10pt }}\cr
e\sb i & 0, & 0, 							& 0, 		& 0, & 0, & 0, &    1,     & 0 &&,\cr
\multispan{11} \hfill\phantom{\vrule height 12pt}
($1$ is at the $(23-i)$-th place for $i=1, \dots, 17$)\hfill  
\cr
\multispan{11}{\dots\dots\dots\phantom{\vrule height 10pt depth 10pt}}\cr
e\sb{18} & 1, & (p+1)/4, & 0, 		& 0, & 0, & 0, &										&	0 &&,\cr
e\sb{19} & 0, & -p,      & 1, 		& 0, & 0, & 0, &										&	0 &&,\cr
e\sb{20} & 0, & 0,       & 0, 		& 1, & -1, & 0, &										&	0 &&,\cr
e\sb{21} & 0, & p,       & 0, 		& 0, & -2, & 0, &										&	0 &&.\cr
}}}
$$
It is easy to see that $h\sb 0\sp 2=-2$,
and that  $e\sb1, \dots, e\sb{21}$ form a basis of $h\sb 0\sperp$.
Since $-(p+1)/2$ is an even integer, 
the $\Z$-condition in this case is as follows:
\begin{equation}\label{eq:ZIV}
x\sb i\in \Z\;\; (i\ne 19, 21), \quad   x\sb {19}, x\sb{21}\in \Z+1/2.
\end{equation}
We use the notation $Q\sb{\sigma}$ 
instead of $Q$ in order to distinguish the situations for  different  Artin invariants $\sigma$.
Suppose that there exist
$a\sb 1, \dots, a\sb{21} \in (1/2)\Z$ satisfying~\eqref{eq:ZIV}
and the inequality  $Q\sb{\sigma}(a\sb 1, \dots, a\sb {21})\le 0$.
We put
$$
J\sb{\sigma,1}:=J\pi (Q\sb{\sigma} (X\sb 1, \dots, X\sb{21}), \{ X\sb{1}\}),
$$
and for $\nu=2, \dots, 17$,
we put
$$
J\sb{\sigma,\nu}:=J\pi (Q\sb{\sigma} (0,\dots, 0, X\sb{\nu} \dots, X\sb{21}), \{ X\sb\nu\}).
$$
Each interval $J\sb{\sigma, \nu}$ is expressed as 
$[-\sqrt{\tau\sb{\sigma, \nu}}, \sqrt{\tau\sb{\sigma, \nu}}]$,
where $\tau\sb{\sigma, \nu}$ are calculated as in Table~\ref{table:tau}.
{\small
\begin{table}
\def\dfrac#1#2{\displaystyle\frac{#1}{#2}}
\renewcommand{\arraystretch}{2.1}
\begin{tabular}{|c||c|c|c|c|c|}
\hline
$\nu\;\;\big\backslash\;\;\sigma$ & $1$ & $3$ & $5$ & $7$ & $9$ \\
\hline
\hline
$1$ & $\dfrac{9p+1}{4p}$ & $\dfrac{7p+3}{4p}$ & $\dfrac{5p+5}{4p}$ & $\dfrac{3p+7}{4p}$ & $\dfrac{p+9}{4p}$ \\[4pt]
\hline
$2$ & $\dfrac{8p+1}{9p+1}$ & $\dfrac{6p+3}{7p+3}$ & $\dfrac{4p+5}{5p+5}$ & $\dfrac{2p+7}{3p+7}$ & $\dfrac{9}{p+9}$ \\[4pt]
\hline
$3$ & $\dfrac{23p+3}{32p+4}$ & $\dfrac{17p+9}{24p+12}$ & $\dfrac{11p+15}{16p+20}$ & $\dfrac{5p+21}{8p+28}$ & $\dfrac{9p+17}{36p}$ \\[4pt]
\hline
$4$ & $\dfrac{14p +2}{23p+3}$ & $\dfrac{10p+6}{17p+9}$ & $\dfrac{6p+10}{11p+15}$ & $\dfrac{2p+14}{5p+21}$ & $\dfrac{8p+8}{9p^2+17p}$ \\[4pt]
\hline
$5$ & $\dfrac{33p +5}{56p+8}$ & $\dfrac{23p+15}{40p+24}$ & $\dfrac{13p+25}{24p+40}$ & $\dfrac{3p+35}{8p+56}$ & $\dfrac{23p+15}{32p^2+32p}$\\[4pt]
\hline
$6$ & $\dfrac{18p +3}{33p+5}$ & $\dfrac{12p+9}{23p+15}$ & $\dfrac{6p+15}{13p+25}$ & $\dfrac{21}{3p+35}$ & $\dfrac{14p+7}{23p^2+15p}$\\[4pt]
\hline
$7$ & $\dfrac{39p+7}{72p+12}$ & $\dfrac{25p+21}{48p+36}$ & $\dfrac{11p+35}{24p+60}$ & $\dfrac{7p+39}{84p}$ & $\dfrac{33p+13}{56p^2+28p}$\\[4pt]
\hline
$8$ & $\dfrac{20p+4}{39p+7}$ & $\dfrac{12p+12}{25p+21}$ & $\dfrac{4p+20}{11p+35}$ & $\dfrac{6p+18}{7p^2+39p}$ & $\dfrac{18p+6}{33p^2+13p}$\\[4pt]
\hline
$9$ & $\dfrac{41p+9}{80p+16}$ & $\dfrac{23p+27}{48p+48}$ & $\dfrac{5p+45}{16p+80}$ & $\dfrac{17p+33}{24p^2+72p}$ & $\dfrac{39p+11}{72p^2+24p}$\\[4pt]
\hline
$10$ & $\dfrac{20p+5}{41p+9}$ & $\dfrac{10p+15}{23p+27}$ & $\dfrac{5}{p+9}$ & $\dfrac{10p+15}{17p^2+33p}$ & $\dfrac{20p+5}{39p^2+11p}$\\[4pt]
\hline
$11$ & $\dfrac{39p+11}{80p+20}$ & $\dfrac{17p+33}{40p+60}$ & $\dfrac{p+9}{20p}$ & $\dfrac{23p+27}{40p^2+60p}$ & $\dfrac{41p+9}{80p^2+20p}$\\[4pt]
\hline
$12$ & $\dfrac{18p+6}{39p+11}$ & $\dfrac{6p+18}{17p+33}$ & $\dfrac{4p+20}{5p^2+45 p}$ & $\dfrac{12p+12}{23p^2+27p}$ & $\dfrac{20p+4}{41p^2+9p}$\\[4pt]
\hline
$13$ & $\dfrac{33p+13}{72p+24}$ & $\dfrac{7p+39}{24p+72}$ & $\dfrac{11p+35}{16p^2+80 p}$ & $\dfrac{25p+21}{48p^2+48p}$ & $\dfrac{39p+7}{80p^2+16p}$\\[4pt]
\hline
$14$ & $\dfrac{14p+7}{33p+13}$ & $\dfrac{21}{7p+39}$ & $\dfrac{6p+15}{11p^2+35p}$ & $\dfrac{12p+9}{25p^2+21p}$ &$\dfrac{18p+3}{39p^2+7p}$\\[4pt]
\hline
$15$ & $\dfrac{23p+15}{56p+28}$ & $\dfrac{3p+35}{84 p}$ & $\dfrac{13p+25}{24p^2+60p}$ & $\dfrac{23p+15}{48p^2+36p}$ &$\dfrac{33p+5}{72p^2+12p}$\\[4pt]
\hline
$16$ & $\dfrac{8p+8}{23p+15}$ & $\dfrac{2p+14}{3p^2+35 p}$ & $\dfrac{6p+10}{13p^2+25p}$ & $\dfrac{10p+6}{23p^2+15p}$ &$\dfrac{14p+2}{33p^2+5p}$\\[4pt]
\hline
$17$ & $\dfrac{9p+17}{32p+32}$ & $\dfrac{5p+21}{8p^2+56 p}$ & $\dfrac{11p+15}{24p^2+40p}$ & $\dfrac{17p+9}{40p^2+24p}$ &$\dfrac{23p+3}{56p^2+8p}$\\[4pt]
\hline
\hline
$20$ & $\dfrac{p+2}{9p+17}$ & $\dfrac{3}{5p+21}$ & $\dfrac{3}{11p+15}$ & $\dfrac{3}{17p+9}$ &$\dfrac{3}{23p+3}$\\[4pt]
\hline
$21$ & $\dfrac{1}{4p+8}$ & $\dfrac{1}{12p}$ & $\dfrac{1}{12p}$ & $\dfrac{1}{12p}$ &$\dfrac{1}{12p}$\\[4pt]
\hline
\end{tabular}
\caption{The table of $\tau\sb{\sigma, \nu}$}
\label{table:tau}
\end{table}
}
We have
\begin{eqnarray}
\label{eq:Jcap1}
J\sb{\sigma, 1}\cap \Z &=& \begin{cases}
\{-1,0,1\} &\textrm{if $\sigma <7$ or $( p, \sigma)=(3, 7), (7, 7)$ or $(3, 9)$},\\ 
\{0\} &\textrm{otherwise};
\end{cases}
\\
\label{eq:Jcap2}
J\sb{\sigma, \nu}\cap \Z &=& \{0\} \quad\textrm{for}\quad \nu=2, \dots, 17.
\end{eqnarray}
From \eqref{eq:Jcap2},
we see that,
if $a\sb 1=0$, then $a\sb 2=\dots=a\sb{17}=0$ holds inductively.
\par
Suppose that $a\sb 1=0$,
so that $a\sb 2=\dots=a\sb{17}=0$.
We put
$$
Q\sprime\sb{\sigma}(X\sb{18}, X\sb{19}, X\sb{20}, X\sb{21}):=
Q\sb{\sigma} (0, \dots, 0, X\sb{18}, X\sb{19}, X\sb{20}, X\sb{21}).
$$
The projection of the $4$-dimensional quadratic body associated with 
$Q\sprime\sb{\sigma}$ to
the $X\sb{20}$-axis is the interval
$[-\sqrt{\tau\sb{\sigma, 20}}, \sqrt{\tau\sb{\sigma, 20}}]$,
where $\tau\sb{\sigma, 20}$ are  given  in Table~\ref{table:tau}.
Hence $a\sb{20}=0$.
The projection of the $3$-dimensional quadratic body associated with 
$Q\sprime\sb{\sigma}(X\sb{18}, X\sb{19}, 0, X\sb{21})$ to
the $X\sb{21}$-axis is the interval
$[-\sqrt{\tau\sb{\sigma, 21}}, \sqrt{\tau\sb{\sigma, 21}}]$,
which is disjoint from $\Z+1/2$.
Thus we get a contradiction to $a\sb{21}\in \Z+1/2$.
\par
From~\eqref{eq:Jcap1},
we have completed the proof except for the cases
where
$\sigma<7$ or $( p, \sigma)=(3, 7), (7, 7), (3, 9)$.
The cases where $p=3$ or $7$ can be treated by  numerical
calculations described in \S\ref{sec:comp}.
Therefore we  assume $p>7$ from now on,
and prove the remaining cases $\sigma<7$.
\par
We have $a\sb 1=\pm 1$.
Suppose that $a\sb 1=1$.
For $\nu\ge 2$, 
we put
$$
[\rho\sp{(\sigma)}\sb{\nu, -}, \rho\sp{(\sigma)}\sb{\nu, +}]
:=
J\pi (Q\sb \sigma  (1,2,0,2,0,\dots, 1+(-1)\sp{\nu-1}, X\sb{\nu}, \dots, X\sb{21}), \{X\sb{\nu}\}).
$$
\par
\smallskip
{\bf Case IV-1}. $\sigma=1$.
The values $\rho\sp{(1)}\sb{\nu, \pm}$ are calculated as in Table~\ref{table:rho}.
Inductively,
we obtain
$$
a\sb{\nu}=1+(-1)\sp{\nu}\qquad\textrm{for $\nu=2, \dots, 10$.}
$$
Since $[\rho\sp{(1)}\sb{11, -}, \rho\sp{(1)}\sb{11, +}]\cap \Z=\emptyset$,
there are no possible values  for $a\sb{11}$.
\par
\smallskip
{\bf Case IV-3}. $\sigma=3$.
Again 
we obtain
$a\sb{\nu}=1+(-1)\sp{\nu}$ for $\nu=2, \dots, 6$
from Table~\ref{table:rho}.
Since $[\rho\sp{(3)}\sb{7, -}, \rho\sp{(3)}\sb{7, +}]\cap \Z=\emptyset$
for $p>7$,
there are no possible  values  for $a\sb{7}$.
\par
\smallskip
{\bf Case IV-5}. $\sigma=5$.
Since
%
%
$$
[\rho\sp{(5)}\sb{2, -}, \rho\sp{(5)}\sb{2, +}]=
\Biggl[
{\dfrac {8\,p+10-\sqrt {4\,{p}^{2}+25\,p+25}}{5\,p+5}}, {\dfrac {8\,p+10+\sqrt
{4\,{p}^{2}+25\,p+25}}{5\,p+5}}
\Biggr],
$$
%
%
we have $a\sb 2=2$.
Since
$$
[\rho\sp{(5)}\sb{3, -}, \rho\sp{(5)}\sb{3, +}]=
\Biggl[
{\dfrac {2\,p-\sqrt {55\,p+75}}{8\,p+10}}, {\dfrac {2\,p+\sqrt {55\,p+75}}{8\,p+10}}
\Biggr],
$$
there are no possible  values  for $a\sb{3}$
except for the cases $p\le 11$.
The case $p=11$ can be treated by numerical calculations.
{\small
\begin{table}
\renewcommand{\arraystretch}{2.3}
\def\dfrac#1#2{\displaystyle\frac{#1}{#2}}
\begin{tabular}{|c|c|c|}
\hline
$\nu$& $\rho\sp{(1)}\sb{\nu,\pm}$ & $\rho\sp{(3)}\sb{\nu,\pm}$ \\ \hline\hline 
$2$& ${\dfrac {16\,p+2\;\pm\;\sqrt {40\,{p}^{2}+13\,p+1}}{9\,p+1}}$ &
${\dfrac {12\,p+6\;\pm\;3\,\sqrt {2\,{p}^{2}+3\,p+1}}{7\,p+3}}$ \\[4pt]\hline
$3$& ${\dfrac {2\,p\;\pm\;\sqrt {92\,{p}^{2}+35\,p+3}}{16\,p+2}}$ &
${\dfrac {2\,p\;\pm\;\sqrt {34\,{p}^{2}+69\,p+27}}{12\,p+6}}$ \\[4pt]\hline
$4$& ${\dfrac {42\,p+6\;\pm\;\sqrt {154\,{p}^{2}+64\,p+6}}{23\,p+3}}$ &
${\dfrac {30\,p+18\;\pm\;\sqrt {50\,{p}^{2}+120\,p+54}}{17\,p+9}}$ \\[4pt]\hline
$5$& ${\dfrac {4\,p\;\pm\;\sqrt {198\,{p}^{2}+96\,p+10}}{28\,p+4}}$ &
${\dfrac {4\,p\;\pm\;\sqrt {46\,{p}^{2}+168\,p+90}}{20\,p+12}}$ \\[4pt]\hline
$6$& ${\dfrac {60\,p+10\;\pm\;\sqrt {234\,{p}^{2}+129\,p+15}}{33\,p+5}}$ &
${\dfrac {40\,p+30\;\pm\;3\,\sqrt {4\,{p}^{2}+23\,p+15}}{23\,p+15}}$ \\[4pt]\hline
$7$& ${\dfrac {6\,p\;\pm\;\sqrt {234\,{p}^{2}+159\,p+21}}{36\,p+6}}$ &
${\dfrac {2\,p\;\pm\;\sqrt {25\,p+21}}{8\,p+6}}$\\[4pt]\hline
$8$& ${\dfrac {70\,p+14\;\pm\;2\,\sqrt {55\,{p}^{2}+46\,p+7}}{39\,p+7}}$ &
\\[4pt]\hline
$9$& ${\dfrac {4\,p\;\pm\;\sqrt {41\,{p}^{2}+50\,p+9}}{20\,p+4}}$ &
\\[4pt]\hline
$10$& ${\dfrac {72\,p+18\;\pm\;\sqrt {100\,{p}^{2}+205\,p+45}}{41\,p+9}}$ &
\\[4pt]\hline
$11$& ${\dfrac {10\,p\;\pm\;\sqrt {195\,p+55}}{40\,p+10}}$ &
\\[4pt]\hline
\end{tabular}
\caption{The table of $\rho\sp{(1)}\sb{\nu,\pm}$ and $\rho\sp{(3)}\sb{\nu,\pm}$}
\label{table:rho}
\end{table}
}
\par
The case where $a\sb{1}=-1$ can be dealt with  in the same way.
\qed
\par
\medskip
Numerical experiments
show that, if  $3< p<10000$,
then we have
$$
R\sb{(X, L)}=\Sigma (h\sb 0\sperp)=A\sb 1 + D\sb{20-2\sigma}.
$$
When $p=3$, we have
$$
R\sb{(X, L)}=\Sigma (h\sb 0\sperp)=A\sb 2 + D\sb{20-2\sigma}.
$$
In the case  $\sigma=1$,
the branch curve $B\sb{(X, L)}$ yields another 
example of a sextic curve with only rational double points
such that  the  total Milnor number is $20$.

\bibliographystyle{amsplain}

\providecommand{\bysame}{\leavevmode\hbox to3em{\hrulefill}\thinspace}

\end{document}